\documentclass[11pt]{article}

\RequirePackage{microtype,enumitem}
\RequirePackage{totpages}
\usepackage{float,amsmath,amssymb,amsthm,bm,balance}
\usepackage{authblk}

\newtheorem{theorem}{Theorem}

\usepackage{algorithm}
\usepackage[noend]{algpseudocode}

\makeatletter
\let\original@algocf@latexcaption\algocf@latexcaption
\long\def\algocf@latexcaption#1[#2]{%
  \@ifundefined{NR@gettitle}{%
    \def\@currentlabelname{#2}%
  }{%
    \NR@gettitle{#2}%
  }%
  \original@algocf@latexcaption{#1}[{#2}]%
}
\makeatother

\usepackage[capitalize]{cleveref}

\def\A{\mathbb{A}}

\def\K{\mathbb{K}}

\def\C{\mathbb{C}}

\def\Q{\mathbb{Q}}
\def\Z{\mathbb{Z}}
\def\N{\mathbb{N}}

\def\softO{O\tilde{~}}

\def\F{\mathbb{F}}

\newcommand{\cS}{{\mathcal S}}

\newtheorem{definition}{Definition}
\newtheorem{proposition}{Proposition}
\newtheorem{lemma}{Lemma}

\AtBeginDocument{%
  }

\begin{document}
\title{$p$-adic algorithm for bivariate Gr\"obner bases}

\author[1]{\'Eric Schost}
\author[2,1]{Catherine St-Pierre}
\date{}

\affil[1]{University of Waterloo, Cheriton School of Computer Science,Waterloo Ontario,Canada}

\affil[2]{Inria (MATHEXP), University Paris-Saclay, Palaiseau, France}

\maketitle


\begin{abstract}
We present a $p$-adic algorithm to recover the lexicographic Gr\"obner
basis $\mathcal G$ of an ideal in $\Q[x,y]$ with a generating set in
$\Z[x,y]$, with a complexity that is less than cubic in terms of the dimension of
$\Q[x,y]/\langle \mathcal G \rangle$ and softly linear in the height
of its coefficients.
We observe that previous results of Lazard's that use Hermite normal
forms to compute Gr\"obner bases for ideals with two generators can be
generalized to a set of $t\in \N^+$ generators. We use this result to obtain a
bound on the height of the coefficients of $\mathcal G$, and to
control the probability of choosing a \textit{good} prime $p$ to build
the $p$-adic expansion of $\mathcal G$.
\end{abstract}



\vspace{-0.5em}
\section{Introduction}

There exists a rich literature dedicated to the solution of polynomial systems in two
variables~\cite{GoKa96,EmTs05,DiEmTs09,AlMoWi10,Rouillier10,BeEmSa11,EmSa12,BoLaPoRo13,LMS13,BoLaMoPoRo14,KoSa14,MeSc16,KoSa15,BoLaMoPoRoSa16,DiDiRoRoSa18,Dahan22},
due in part to their numerous applications in real algebraic geometry
and computer-aided design. Our focus in this paper is on the
complexity of computing the lexicographic Gr\"obner basis of a
zero-dimensional ideal in $\Q[x,y]$, specifically by means of $p$-adic
techniques based on Newton iteration. An important aspect of this work
is to give bit-size bounds for such a Gr\"obner basis, as well as bounds
on the number of primes of bad reduction.

$p$-adic techniques have been considered in the context of Gr\"obner
basis computations (in an arbitrary number of variables) for
decades. In 1983 and 1984, Ebert and Trinks addressed the question of
modular algorithms for Gr\"obner bases~\cite{Ebert83,Trinks84},
specifically for systems without multiple roots; these techniques were also
used in geometric resolution
algorithms~\cite{GiHeMoPa95,GiHeMoMoPa98,GiHaHeMoMoPa97,GiLeSa01}. The
absence of multiple roots allows for simple and efficient algorithms;
for arbitrary inputs, the question is more involved.

Winkler gave the first $p$-adic algorithm to construct a Gr\"obner
basis~\cite{Winkler88} that applies to general inputs; Pauer refined
the discussion of good prime numbers~\cite{Pauer92}, and Arnold
revisited, and simplified, these previous constructions
in~\cite{Arnold03}. No complexity analysis was provided; these
$p$-adic algorithms remain complex (they not only lift the Gr\"obner
basis, but also the transformation matrix that turns the input system
into its Gr\"obner basis), and to our knowledge, achieve linear
convergence.

In the specific context of bivariate equations, $p$-adic techniques
have already been put to use in previous work, first in the particular
case of non-multiple roots~\cite{LMS13}, then to compute a
set-theoretic description of all roots, even in the presence of
multiplicities~\cite{MeSc16}. However, in this case, the latter
algorithm does not reveal the local structure at multiple roots.

In~\cite{ScSt23}, we presented a form of Newton iteration specifically
tailored to lexicographic Gr\"obner bases in two variables. It crucially rests on results due to Conca
and Valla~\cite{CoVa08}, who gave an explicit parametrization of
bivariate ideals with a given initial ideal: our lifting algorithm
works specifically with the parameters introduced by Conca and
Valla. Our contribution in this paper is to build on~\cite{ScSt23} to
give a complete $p$-adic algorithm: we quantify bad primes, show how
to initialize the lifting process, give bounds on the size of the
output, and analyze the cost of the whole algorithm.

The following theorem gives a slightly simplified form of our main
result, where the probability of success and the number of input
polynomials are kept constant (the more precise version is given in
the last section). In what follows, the {\em height} of a nonzero
integer $u$ is $\log(|u|)$; if $\mathcal G$ is a family of polynomials
in $\Q[x,y]$, we define $\deg(\mathcal G)=\dim_\Q \Q[x,y]/\mathcal{G}$
and let $h(\mathcal G)$ be the maximum height of the numerators and
denominators of its coefficients.

\begin{theorem}\label{cor:Zcomplexity}
  Let $\mathcal F=(f_1,\dots,f_t)$ be in $\Z[x,y]$, with degree at
  most $d$, height at most $h$, and with finitely many common
  solutions in $\C^2$. Let $\mathcal G$ be the lexicographic Gr\"obner
  basis of $I$ for the order $x \prec y$ and write
  $s=|\mathcal G|$, $\delta=\deg(\mathcal G)$, $b=h(\mathcal G).$

  For $P>0$, assuming $P \in O(1)$ and $t \in O(1)$, there is an
  algorithm that computes $\mathcal G$ with probability of success at
  least $1-1/2^P$ using a number of bit operations softly linear in
  $$d^2 h + ( d^{\omega+1} + \delta^\omega)\log(h) + (d^2 \delta + d
  \delta^2 + s^2\delta^2) (b+\log(h)).$$
\end{theorem}
With the notation in the theorem, the bitsize of the input
is linear in $d^2h$, and that of the output is linear in $s\delta b$.
If all solutions of $\mathcal F$ have multiplicity one, previous forms
of Newton iteration achieve better runtimes, softly linear in the
output size~\cite{GiLeSa01,Schost03,LMS13,DMSWX05}
 (but instead of a
Gr\"obner basis, they compute a {\em triangular decomposition} of
$V(\mathcal F)$, or change coordinates). Hence, it makes sense to apply our techniques only to multiple solutions.

This is what provides the motivation for our second result, where we
compute the Gr\"obner basis of the $\langle x,y\rangle$-primary
component $J$ of $I$. As a natural extension, one can consider
combining this with~\cite{HyMeScSt19}, which shows how to put an
arbitrary primary component of $I$ in correspondence with the $\langle
x,y\rangle$-primary component of a related ideal in $\K[x,y]$ for a
finite extension $\K$ of $\Q$.

\begin{theorem}\label{cor:Zcomplexity0}
  Let $\mathcal F=(f_1,\dots,f_t)$ be in $\Z[x,y]$, with degree at
  most $d$, height at most $h$, with finitely many common
  solutions in $\C^2$. Let $\mathcal G^0$ be the lexicographic
  Gr\"obner basis of the $\langle x,y\rangle$-primary component of $I$
  for the order $x \prec y$ and write
  $r=|\mathcal G^0|$, $\eta=\deg(\mathcal G^0)$,  $c=h(\mathcal G^0).$
    For $P>0$, assuming $P \in O(1)$ and $t \in O(1)$, there is an
  algorithm that computes $\mathcal G^0$ with probability of success at
  least $1-1/2^P$ using a number of bit operations softly linear in
  $$d^2 h + (d^\omega\eta+\eta^\omega)\log(h) + \eta^2 c.$$
\end{theorem}

\noindent{\bf Outlook.} Inspired by Lazard \cite{Lazard85}, we prove
in \cref{coeff} that the Hermite Normal form of an ``extended
Sylvester matrix'' built from $f_1,\dots,f_t$ gives the coefficients
of a {\em detaching basis} of the ideal $I$ they generate. We also
present a variant of this result, where replacing the Hermite normal form
by the Howell normal form yields a Gr\"obner basis of a localization of
$I$. We use these results in two manners: to compute the initial
Gr\"obner basis modulo $p$, prior to entering Newton iteration, and to
obtain bit-size bounds for the output (over $\Q$) and quantify bad
choices of the prime $p$. The underlying algorithms for the above theorems are in~\cref{sec:mainalgo}.

\vspace{-0.5em}
\section{Using matrix normal forms}\label{coeff}

In this section, we assume $I = \langle f_1,\dots,f_t\rangle \subset
\K[x,y]$, for $t\geq 2$, and we derive the lexicographic Gr\"obner
basis of $I$, or its primary component at the origin, from either
Hermite or Howell normal forms of matrices over $\K[x]$, for an
arbitrary field $\K$. These results are direct extensions of previous
work of Lazard's~\cite{Lazard85}, who used Hermite forms in the case
$t=2$.

In what follows, for a subset $S \subset \K[x,y]$ and $n \ge 0$, we
let $S_{<(.,n)}$ be the subset of all $f$ in $S$ with $\deg_y(f) < n$;
notation such as $S_{\le(.,n)}$ is defined similarly.  In particular,
if $S$ is an ideal of $\K[x,y]$, $S_{<(.,n)}$ is a free $\K[x]$-module
of rank at most $n$. For $S=\K[x,y]$ itself, $\K[x,y]_{<(.,n)}$ is a
free $\K[x]$-module of rank $n$, equal to $\bigoplus_{0 \le i < n}
\K[x]y^i$.
For such an $n$, we also let $\pi_n$ denote the $\K[x]$-module isomorphism
$\K[x,y]_{<(.,n)} \to \K[x]^n$, which maps $f_0 + \cdots +
f_{n-1}y^{n-1}$ to the vector $[f_{n-1}\cdots f_0]^\top$.

\vspace{-0.5em}
\subsection{Detaching bases}

Let $I$ be an ideal in $\K[x,y]$ and let $\mathcal G=(g_0,\dots,g_s)$
be its reduced minimal Gr\"obner basis for the lexicographic order
induced by $y\succ x$, listed in decreasing order; we write
$n_i=\deg_y(g_i)$ for all $i$ (so these exponents are decreasing). We
define polynomials $A_0,A_1,\dots$ as follows: for $0 \le i < n_s$,
$A_i=0$, and if there exists $k$ in $\{0,\dots,s\}$ such that $n_k
=i$, $A_i=g_k$; otherwise, $A_i$ is obtained by starting from $y
A_{i-1}$, and reducing all its terms of $y$-degree less than $i$ by
$\mathcal G$.

For example, if $I$ has a Gr\"obner basis of the form $(y-f(x),
g(x))$, the polynomials $A_i$ are given by $A_0=g$, $A_1=y-f$ and for
$i \ge 2$, $A_i= y^i -(f^i \bmod g)$ (see~\cite{Ayoub83}
for a previous discussion).

\begin{lemma}
  For $i \ge n_s$, $\deg_y(A_i)=i$.
\end{lemma}
\vspace{-0.5em}
\begin{proof}
  This is true for $i$ of the form $n_k$. For $i$ in
  $n_{k},\dots,n_{k-1}-1$, we proceed by induction, with the remark
  above establishing the base case (for $k=0$, we consider all $i \ge
  n_0$). Assume $\deg_y(A_{i-1})=i-1$, so that $\deg_y(y A_{i-1})=i$.
  Because we use the lexicographic order $x \prec y$, the reduction of
  the terms of $y$-degree less than $i$ in $yA_{i-1}$ does not
  introduce terms of $y$-degree $i$ or more.
\end{proof}

\vspace{-4mm}
\begin{lemma}\label{lemma:freem}
  For $n \ge n_s$, the $\K[x]$-module $I_{\le (.,n)}$ is free of rank
  $n-n_s+1$, with basis $A_{n_s},\dots,A_n$.
\end{lemma}
\vspace{-5mm}
\begin{proof}
  The polynomials $A_{n_s},\dots,A_n$ are all nonzero, with pairwise
  distinct $y$-degrees, so they are $\K[x]$-linearly independent. They
  all belong to $I_{\le (.,n)}$, so it remains to prove that they
  generate it as a $\K[x]$-module.
  This is done by induction on $n \ge n_s$. Take $f$ in $I_{\le
    (.,n)}$, and write it as $f=f_n y^n + g$, with $f_n$ in $\K[x]$
  and $g$ in $K[x,y]_{\le (.,n-1)}$. Let $h \in \K[x]$ be the
  polynomial coefficient of $y^n$ in $A_n$, so that $A_n = h_n y^n +
  B_n$, with $B_n$ in $\K[x,y]_{\le (.,n-1)}$. Write the Euclidean
  division $f_n = q h_n + r$ in $\K[x]$, with $\deg_x(r) < \deg_x(h_n)$:
  \begin{align*}
    f&=(q h_n +r)y^n  + g=qh_n y^n + r y^n +g= q A_n- q B_n +r y^n + g.
  \end{align*}
  The polynomial $- q B_n +r y^n + g$ is in $I$, so its normal form
  modulo $\mathcal G$ is zero. The terms $- q B_n+g$ have $y$-degree
  less than $n$, so their normal form has $y$-degree less than $n$ as
  well; since $r y^n$ is already reduced modulo $\mathcal G$, it must
  be zero.
  It follows that $f=q A_n +g - q B_n$, with $g-q B_n$ in $I_{\le
    (.,n-1)}$. If $n=n_s$, this latter polynomial must vanish, proving the base case. Else, by induction
  assumption, it is a $\K[x]$-linear combination of
  $A_{n_s},\dots,A_{n-1}$.
\end{proof}
\vspace{-0.5em}
For $n \ge n_0$, the {\em detaching basis} of $I$ in degree $n$ is the
sequence $(A_{n_s},\dots,A_{n})$. Because we take $n \ge n_0$, this is
(in general) a non-minimal Gr\"obner basis of $I$, and we can recover
$\mathcal G$ from it by discarding redundant entries (that is, all
polynomials whose leading term is a multiple of another leading term).

\vspace{-0.5em}
\subsection{Using Hermite normal forms}\label{sec: Hermite}

Given $\mathcal{F}=(f_1,\dots,f_t)$ in $\K[x,y]$, we prove that the
Hermite normal form of a certain Sylvester-like matrix associated to
them gives a lexicographic detaching basis of the ideal $I$ they
generate. In~\cite{Lazard85}, Lazard covered the case $t=2$, under an
assumption on the leading coefficients (in $y$) of the $f_i$'s.

We extend his work (in a direct manner) to situations where such
assumptions do not hold. First, to polynomials $\mathcal F=
(f_1,\dots,f_t)$ in $\K[x,y]$, we associate an integer
$\Delta(\mathcal{F})$, defined as follows.
\vspace{-0.4em}
\begin{definition}
  Let $\mathcal F= (f_1,\dots,f_t)\in \K[x,y]^t$, with
  $(A_{n_s},\dots,A_{n_0})$ as detaching basis in degree $n_0$, with
  $n_0$ and $n_s$ the maximal, resp.\ minimal $y$-degree of the
  polynomials in the lexicographic Gr\"obner basis of $\langle
  f_1,\dots,f_t \rangle$, for the order $x \prec y$.
  
  We let $\Delta(\mathcal F)$ be the minimal integer $\Delta$ such
  that for $i=n_s,\dots,n_0$, there exist $w_{i,1},\dots,w_{i,t}$ in
  $\K[x,y]^t$, all of $y$-degree less than $\Delta$, and such that
  $A_i = w_{i,1} f_1 + \cdots + w_{i,t} f_t$.
\end{definition}
The following proposition gives the basic result using this definition, allowing us to extract a detaching basis from a Hermite
form computation. We use {\em column} operations, with Hermite normal
forms being lower triangular. The first nonzero entry in a nonzero
column is called its {\em pivot}, its index being called the {\em
  pivot index}. Pivots 
are monic in $x$.
\vspace{-0.5em}
\begin{proposition}\label{detached}
  Let $\mathcal F=(f_1,\dots,f_t)$ be in $\K[x,y]$, for $t \ge 2$, of
  $y$-degree at most $d_y$, and assume that they generate an ideal
  $I=\langle f_1,\dots,f_t\rangle$ of dimension zero. For
  $i=1,\dots,t$, write $f_i = f_{i,0} + \cdots + f_{i,d_y}y^{d_y}$,
  with all $f_{i,j}$ in $\K[x]$.

  For $D \ge \Delta(\mathcal F)$, let $c_1,\dots,c_{K}$ be the nonzero
  columns of the Hermite normal form $\bm H$ of $\bm S = [\bm
    S_1\cdots \bm S_t] \in \K[x]^{(d_y+D) \times tD}$, where $$ \bm
  S_i= \begin{bmatrix} f_{i,d_y} & &\\ \vdots & \ddots\\ f_{i,0} &
    &f_{i,d_y}\\ &\ddots &\vdots\\& &f_{i,0} \\ \end{bmatrix} \in
  \K[x]^{(d_y+D)\times D}.$$ Then, there exists $K' \le K$ such that
  $\pi_{d_y+D}^{-1}(c_{K'})$ is monic in $y$; with $K'$ the largest
  such integer, $\pi_{d_y+D}^{-1}(c_K),\dots,\pi_{d_y+D}^{-1}(c_{K'})$
  is a detaching basis of $I$.
\end{proposition}

Thus, while we do not know the $y$-degrees $n_i$ of the
elements in the Gr\"obner basis of $I$, as long as $D \ge
\Delta(\mathcal F)$, it is enough to consider the last nonzero columns
of $\bm H$, stopping when we find (through $\pi_{d_y+D}^{-1}$) a
polynomial that is monic in $y$. 
\vspace{-2mm}
\begin{proof}
  Let $D \ge \Delta(\mathcal F)$ be as in the proposition. Let us index
  the columns of each block $\bm S_i$ by $y^{D-1},\dots,y,1$, and its
  rows by $y^{d_y+D-1},\dots,y,1$. Then, $\bm S_i$ is the matrix of the
  map $\K[x,y]_{< (.,D)} \to \K[x,y]_{< (.,d_y+D)}$ given by $w_i
  \mapsto w_i f_i$. The matrix $\bm S$ itself maps a vector
  $(w_1,\dots,w_t)$, with all entries of $y$-degree less than $D$, to
  $\sum_{i=1}^t w_i f_i \in I_{<(.,d_y+D)}$.

  Let $\mathcal G=(g_0,\dots,g_s)$ be the lexicographic Gr\"obner basis of
  $I=\langle f_1,\dots,f_t\rangle$, $f_i\succ f_{i+1}$, with
  $\deg_y(g_i) = n_i$ for all $i$. Since we assume that $I$ has
  dimension zero, we have $n_s=0$, and $g_0$ is monic in $y$.
  
  Let $A_0,\dots,A_{n_0}$ be the detaching basis of $I$ in degree
  $n_0$. We denote by $c_1,\dots,c_K$ the nonzero columns of the
  Hermite form $\bm H$ of $\bm S$, and we let
  $H_i=\pi_{d_y+D}^{-1}(c_{i})$, for $i=0,\dots,n_0$. We will prove that
  $A_{i} = H_{K-i}$ for $i=0,\dots,n_0$. Since $g_0$ is the only
  polynomial in $A_0,\dots,A_{n_0}$ which is monic in $y$, this
  will
  establish the proposition, with $K' = K - n_0$.
  
  Since both $A_i$ and $H_{K-i}$ are in $I$, to prove that they are
  equal, it is enough to prove that for all $i$, $A_i-H_{K-i}$ is
  reduced with respect to the Gr\"obner basis $\mathcal G$ of $I$.
  Because $D \ge \Delta(\mathcal F)$, we deduce that $A_0,\dots,A_{n_0}$
  are in the column span of $\bm S$. Since they have respective
  $y$-degrees $0,\dots,n_0$, we see that $\deg_y(H_{K-i}) =
  \deg_y(A_i)=i$ for all $i=0,\dots,n_0$. In addition, for all such
  $i$, we can write $A_i = \sum_{j=0}^{i} a_{i,j} H_{K-j}$, for some
  $a_{i,j}$ in $\K[x]$.

  However, \cref{lemma:freem} shows that for the same index
  $i$, we can write $H_{K-i} = \sum_{j=0}^i b_{i,j} A_j$, for some
  $b_{i,j}$ in $\K[x]$. Because both $A_i$ and $H_{K-i}$ have leading
  $y$-coefficients that are monic in $x$, it follows that
  $b_{i,i}=a_{i,i}=1$ for all $i$. This proves that $A_i$ and
  $H_{K-i}$ have the same coefficient of $y$-degree $i$ (call it $M_i
  \in \K[x]$), and thus that $A_i-H_{K-i}$ has $y$-degree less than
  $i$.

  By the definition of a detaching basis, all terms of $y$-degree less
  than $i$ in $A_i$ are reduced with respect to $\mathcal G$. On the other
  hand, by the property of Hermite forms, for $j < i$, the coefficient
  of $y$-degree $j$ in $H_{K-i}$ is reduced with respect to
  $M_j$. Since we saw that $M_j$ is also the coefficient of $y^j$ in
  $A_j$, this proves that all terms of $y$-degree less than $i$ in
  $H_{K-i}$ are reduced with respect to $A_0,\dots,A_{i-1}$, and thus
  with respect to $\mathcal G$. Altogether, $A_i-H_{K-i}$ itself is reduced
  with respect to $\mathcal G$, which is what we set out to prove.
\end{proof}

We call $\textsc{HermiteGroebnerBasis}(\mathcal F, D)$ a procedure
that takes as input $\mathcal F=(f_1,\dots,f_t)$ and $D$, and returns
the lexicographic Gr\"obner basis of $I=\langle f_1,\dots,f_t \rangle$
obtained by computing the Hermite normal form of $\bm S$ as above,
extracting the Gr\"obner basis of $I$ from its detaching basis. Here,
we take for $d_y$ the maximum degree of the $f_i$'s, and we assume
that we have $D \ge \Delta(\mathcal F)$ and $D \ge d_y$.

The assumption that the ideal $I$ has dimension zero implies that it
contains a non-zero polynomial in $\K[x]$; as a result, its detaching
basis has entries of $y$-degrees $0,1,\dots$, so that the Hermite form
of ${\bm S}$ is lower triangular with $d_y+D$ non-zero diagonal
entries. In other words, $\bm S$ has rank $d_y+D$ (seen as a matrix
over $\K(x)$).

If $t=2$ and $D=d_y$, this matrix is square, but in general, it may
have more columns than rows (recall that we assume $D \ge d_y$). Using
the algorithm of~\cite{LaNeVuZh22}, we can permute the columns of $\bm
S$ to find a $(d_y+D) \times tD$ matrix $\bm S'$ whose leading
$(d_y+D)\times (d_y+D)$ minor is nonzero; this takes $\softO(tD^\omega
d)$ operations in $\K$, with $d$ the maximum degree of the $f_i$'s. Let
us define the $tD \times tD$ square matrix
\begin{equation}\label{eq:defSsq}
{\bm S}^{{\rm sq}}
= \begin{bmatrix}{\bm S}'\\
\begin{matrix}
\bm 0_{(t-1)D-d_y,d_y+D} &\bm I_{(t-1)D-d_y,(t-1)D-d_y}
\end{matrix}\end{bmatrix}
\end{equation}
together with its Hermite form ${\bm H}^{{\rm sq}}$; the first $d_y+D$
rows of it give us the Hermite form $\bm H$ of~$\bm S$. The Hermite
form of ${\bm S}^{{\rm sq}}$ is computed in $\softO(t^\omega
D^{\omega} d)$ operations in $\K$~\cite{LaNeZh17}. This gives the
overall cost of computing the lexicographic Gr\"obner basis of $I$,
assuming an upper bound on $\Delta(\mathcal F)$ is known. To our
knowledge, not much exists in the literature on complete cost analysis
for bivariate ideals, apart from Buchberger's algorithm, with cost
$\frac{3}{2}(t+2(d+2)^2)^4$~\cite{buchberger1983note}.

The following proposition gives various bounds on $\Delta(\mathcal
F)$, whose strength depends on the assumptions we make on $\mathcal
F$. The first one is a direct extension of
Lazard's~\cite[Lemma~7]{Lazard85}, and is linear in the $y$-degree of
the input. The others are based on results
from~\cite{Kollar88,DiFiGiSe91}, which involve total degree
considerations.
\vspace{-0.5em}
\begin{proposition}\label{deg span}\label{DeltaNullst}
  Let $\mathcal F= (f_1,\dots,f_t)$ be in $\K[x,y]$ of degree at most
  $d \ge 1$, and $y$-degree at most $d_y$, and let $I=\langle
  f_1,\dots,f_t\rangle \subset \K[x,y]$. Set $d'=\max(d,3)$. Then:
  \begin{itemize}[leftmargin=*]
  \item if there exists $i$ in $\{1,\dots,t\}$ such that the
    coefficient of $y^d$ in $f_i$ is a nonzero constant, $\Delta(\mathcal
    F) \le \Delta_1(d_y)=d_y$
  \item if $t=2$ and $I$ has finitely many zeros over $\overline \K$,
    $\Delta(\mathcal F) \le \Delta_2(d) = 2{d'}^2+d' \in O(d^2)$
  \item if $I$ has finitely many zeros over $\overline \K$,
    $\Delta(\mathcal F) \le \Delta_3(d)= 16{d'}^4 +2{d'}^2+2d' \in O(d^4)$
  \end{itemize}
\end{proposition}
\vspace{-4mm}
\begin{proof}[First item.]
  In what follows, without loss of generality, we assume that the
  coefficient of $y^{d_y}$ in $f_t$ is $1$. We prove a slightly more
  general claim: {\em any polynomial} $f$ in $I_{<(.,2d_y)}$ can be
  written as $f= w_1 f_1 + \cdots + w_t f_t$, with all $w_i$ in
  $\K[x,y]_{<(.,d_y)}$. This is enough to conclude since all entries
  $A_{n_s},\dots,A_{n_0}$ in the detaching basis of $I$ in degree
  $n_0$ have $y$-degree at most $d_y \le 2d_y-1$ (because we
  use a lexicographic order with $x \prec y$).

  Let thus $f$ be given in $I_{<(.,2d_y)}$.  There is a family
  $w=(w_1,\dots,w_t)$ in $\K[x,y]$ such that $f =
  \sum_{i=1}^tw_{i}f_i$, since $f$ is in $I$. For such a family $w$,
  we define $\cS_w = \{i\mid \deg_y (w_i) \ge d_y\}$.  For any $w$
  such that $\cS_w$ is not empty, we further set $\nu_w= \min(\cS_w)
  \in \{1,\dots,t\}$, and we let $\nu$ be the {\em maximal} value of
  these $\nu_w$'s. It is well-defined, since there is a vector $w$ for
  which $\cS_w$ is not empty (we can replace $(w_{t-1},w_t)$ by
  $(w_{t-1}+g f_t,w_t-gf_{t-1})$ for any $g$ in $\K[x,y]$).
  
  Let $w$ be such that $\nu=\nu_w$. We claim that $\cS_w \neq \{t\}$:
  otherwise we would have $\deg_y (w_tf_t) \ge 2d_y$, while
  $\deg_y(w_if_i) < 2d_y$ for all other $i$'s; this would contradict
  the assumption $\deg_y(f) < 2d_y$. This shows that $\nu < t$.
  Let us further refine our choice of $w$, by taking it such that,
  among all those vectors for which $\cS_w$ is not empty and
  $\nu_w=\nu$, the $y$-degree of $w_\nu$ is minimal. Let us then write
  $e=\deg_y(w_\nu)$ (so that $e \ge d_y$) and let $c \in \K[x]$ be the
  coefficient of $y^e$ in $w_\nu$. We can use it to rewrite $f$ as
  $$ f = \sum_{i=1}^tw_{i}f_i + c y^{e-d_y} f_\nu f_t - c
  y^{e-d_y}f_tf_\nu.$$ If we set 
  $$ w_i'= \begin{cases} w_\nu - c
    y^{e-d_y}f_t&\text{ when } i=\nu;\\ w_t + c y^{e-d_y}f_\nu&\text{ when
    } i=t;\\ w_i&\text{ otherwise,} \\ \end{cases} $$ 
    we still have $
  f = \sum_{i=1}^t w_i'f_i.$ By construction, $\deg_y(w'_i) =
  \deg_y(w_i) < d_y$ for all $i < \nu$, so none of $1,\dots,\nu-1$ is in
  $\cS_{w'}$. If $\nu$ is in $\cS_{w'}$, then the inequality
  $\deg_y(w'_\nu) < \deg_y(w_\nu)$ contradicts the choice of $w$, so
  that $\nu$ is not in $\cS_{w'}$. This shows that $\cS_{w'}$ is empty,
  since otherwise its minimum element would be greater than $\nu$.
\end{proof}

\vspace{-2mm}
For the second and third items, we use results from~\cite{DiFiGiSe91},
for which we need total degree bounds on the input polynomials
$\mathcal F=(f_1,\dots,f_t)$ and the elements $A_0,\dots,A_{n_0}$ in
the detaching basis (here $n_s=0$ since $I$ having finitely many
solutions implies that it contains a nonzero polynomial in $\K[x]$).
For the inputs $f_i$, we have the degree bound $\deg(f_i) \le d \le
d'$.  For the $A_i$'s, we have the bounds $\deg_x(A_i) \le d^2$ (by
B\'ezout's theorem) and $\deg_y(A_i) \le d$ for $i \le n_0$, so their
total degree is at most $D=d'{}^2+d'$.

\vspace{-2mm}
\begin{proof}[Second item.]
  When $t=2$ and $I$ has dimension zero (that is, has a finite,
  nonzero number of solutions in $\overline \K$), $f_1,f_2$ are in
  complete intersection, so that we have $A_i = w_{i,1} f_1 + w_{i,2}
  f_2$, with $\deg_y(w_{i,j}) \le D+{d'}^2$ for all $i,j$, by
  Theorem~5.1 in~\cite{DiFiGiSe91}. Overall, the resulting degree
  bound is $2{d'}^2+ d'$.
    If we assume that $I =\K[x,y]$, we know that there
  are $g_1,g_2$ in $\K[x,y]$ such that $g_1 f_1 + g_2 f_2=1$, with
  $\deg(g_i) \le {d'}^2$~\cite{Kollar88}. Multiplying this by $A_j$,
  for $j \le n_0$, we obtain the expression $(g_1 A_j) f_1 + (g_2 A_j)
  f_2=A_j$, with $\deg_y(g_i A_j) \le {d'}^2+d$ in this case.
\end{proof}

\vspace{-4mm}
\begin{proof}[Third item.] \cite[Corollary~3.4]{DiFiGiSe91} gives equalities $A_i =
  w_{i,1} f_1 + \cdots + w_{i,t} f_t$, with $\deg_y(w_{i,j}) \le
  D+16{d'}^4 +{d'}^2+d'$ for all $i,j$.
\end{proof}
\vspace{-0.5em}
\vspace{-0.5em}
\subsection{Using the Howell form}

We now investigate how using another matrix normal form, the {\em
  Howell} form~\cite{Howell86}, yields information about certain
primary components of an ideal $I$ as above.
Howell forms are defined for matrices with entries in a principal
ideal ring $\A$; below, we will take $\A=\K[x]/x^k$, for an integer
$k$. Again, we consider column operations then an
$n \times m$ matrix $\bm H$ over $\A=\K[x]/x^k$ is in Howell normal
form if the following (taken from~\cite[Chapter~4]{Storjohann00})
hold:
\begin{enumerate}[leftmargin=*]
\item let $r \le m$ be the number of nonzero columns in $\bm H$; then
  these nonzero columns have indices $1,\dots,r$
\item $\bm H$ is in lower echelon form: for $i=1,\dots,r$, let $j_i
  \in \{1,\dots,n\}$ be the index of the first nonzero entry in the
  $i$th column; then, $j_1 < \cdots < j_r$
\item all pivots $H_{j_i,i}$, for $i=1,\dots,r$, are of the form
  $x^{c_i}$
\item for $i=1,\dots,r$ and $k=1,\dots,i-1$, $H_{j_i,k}$ is reduced modulo
  $H_{j_i,i}$
\item for $i=0,\dots,r$, any column in the column span of $\bm H$ with
  at least $j_i$ leading zeros is an $\A$-linear combination of
  columns of indices $i+1,\dots,r$ (here, we set $j_0=0$)
\end{enumerate}

For any $\bm M$ in $\A^{n \times m}$, there is a unique $\bm H$
in Howell normal form in $\A^{n \times m}$, and a not necessarily unique
invertible matrix $\bm U$ in $\A^{m\times m}$ such that $\bm H = \bm M
\bm U$. The matrix $\bm H$ is the Howell normal form of $\bm M$.

Given $f_1,\dots,f_t$ as before, we are interested here in computing
the lexicographic Gr\"obner basis of $J=\langle
f_1,\dots,f_t,x^k\rangle$, for a given integer $k$. In particular, if
$(0,0)$ is in $V(f_1,\dots,f_t)$, and no other point $(0,\beta)$ is,
for $\beta\ne 0$, $J$ is the $\langle x,y\rangle$-primary component of
$I=\langle f_1,\dots,f_t\rangle$, if $k$ is large enough.

The following proposition shows how to reduce this computation to a
Howell normal form calculation. In what follows, the {\em canonical
  lift} of an element in $\A=\K[x]/x^{k}$ to $\K[x]$ is its unique
preimage of degree less than $k$; this carries over to vectors and
matrices (and in particular to the output of the Howell form
computation). Contrary to what happens for Hermite forms, there is no
guarantee that the polynomials extracted from the Howell form are a
detaching basis, as we may be missing the first polynomial (that
belongs to $\K[x]$) and its multiples. The proposition below restores
this by considering a few extra columns, if needed.

\begin{proposition}
  Let $f_1,\dots,f_t$ be in $\K[x,y]$, for $t \ge 2$, of $y$-degree at
  most $d_y$, and assume that they generate an ideal of dimension
  zero. Let $k$ be a positive integer and $\A=\K[x]/x^{k}$.

  For $D \ge \Delta(f_1,\dots,f_t,x^k)$, let $\bm B \in \A^{(d_y+D)
    \times tD}$ be the Howell normal form of $\bar{\bm S} =\bm S \bmod
  x^k$, with $\bm S$ as in Proposition~\ref{detached}, and let $\bm
  B_{\rm lift}$ be its canonical lift to $\K[x]^{(d_y+D)\times tD}$.

  Let $h_1,\dots,h_L$ be the nonzero columns of $\bm B_{\rm lift}$,
  and let $r \in \{1,\dots,d_y+D\}$ be the pivot index of $h_L$. Set
  $L''=L+d_y+D-r$ and, for $i=L+1,\dots,L''$ let
  $h_i=[0~\cdots~0~x^k~0~\cdots~0]^\top$, with $x^k$ at index $r+i-L
  \in\{r+1,\dots,d_y+D\}$.

  Then, there exists $L' \le L$ such that $\pi_{d_y+D}^{-1}(h_{L'})$ is
  monic in $y$; with $L'$ be the largest such integer,
  $\pi_{d_y+D}^{-1}(h_{L''}),\dots,\pi_{d_y+D}^{-1}(h_{L'})$ is 
  a detaching basis of $\langle f_1,\dots,f_t,x^k\rangle$.
\end{proposition}
\vspace{-0.5em}
\begin{proof}
  Let $\Gamma=(\Gamma_0,\dots,\Gamma_\sigma)$ be the lexicographic
  Gr\"obner basis of $J=\langle f_1,\dots,f_t,x^k\rangle$, listed in
  decreasing order, with $\Gamma_i$ of $y$-degree $\nu_i$ for all $i$;
  since $x^k$ is in $J$, $\nu_\sigma=0$. Then, let
  $C_{0},\dots,C_{\nu_0}$ be the detaching basis of $J$ in degree
  $\nu_0$, with $\deg_y(C_i)=i$ for all $i$.

  We know that the first polynomials in the detaching basis are of the
  form $C_0=x^\ell,C_1=y
  x^\ell,\dots,C_{\nu_{\sigma-1}-1}=y^{\nu_{\sigma-1}-1}x^\ell$, for
  some $\ell \le k$. If $\ell =k$, they all vanish modulo $x^k$, but
  the next polynomial $C_{\nu_{\sigma-1}}$ does not. If $\ell < k$,
  none of them vanishes modulo $x^k$. Thus, we define
  $\rho=\nu_{\sigma-1}$ in the former case and $\rho=0$ in the
  latter. 

  Let further $D \ge \Delta(f_1,\dots,f_t,x^k)$ be as in the
  proposition. If we consider the extended Sylvester matrix $\bm T \in
  \K[x]^{(d_y+D) \times (t+1)D}$ built from $f_1,\dots,f_t,x^{k}$, the
  assumption on $D$ shows that each $\pi_{d_y+D}(C_i)$ is in the
  column span of $\bm T$.  For $i=0,\dots,\nu_0$, we let $v_i$ be the
  column vector $\pi_{d_y+D}(C_i) \bmod x^k \in \A^{d_y+D}$; the
  previous paragraph shows that the nonzero vectors $v_i$ are
  precisely $v_\rho,\dots,v_{\nu_0}$. By reduction modulo $x^k$ of the
  membership relations above, we see that $v_\rho,\dots,v_{\nu_0}$ are
  in the $\A$-span of the columns of $\bar{\bm S}$.
  
  Lazard's structure theorem~\cite[Theorem~1]{Lazard85} shows that every polynomial
  $\Gamma_j$ in the reduced Gr\"obner basis of $J$ is of the form
  $\Gamma_j=x^{m_j} \gamma_j$, with $\gamma_j$ monic in $y$ and $m_j
  \le \ell$ (the inequality is strict, except for $j=0$). It follows
  that for $i=\rho,\dots,\nu_0$, the pivot in $v_i$ is also a power of
  $x$, at index $d_y+D-i$ (precisely, it is $x^{m_j}$, for $j$ the
  largest integer such that $\nu_j \le i$).
  
  Let $\eta_1,\dots,\eta_L$ be the nonzero columns in the Howell form
  $\bm B$ of $\bar{\bm S}$. By definition of the Howell form, the
  former observation implies that for $i=\rho,\dots,\nu_0$, $v_i$ is
  in the $\A$-span of those $\eta_j$'s starting with at least
  $d_y+D-i-1$ zeros. For such an $i$, since the entry at index
  $d_y+D-i$ in $v_i$ is nonzero, there  exists (exactly) one
  $\eta_j$ with pivot index $d_y+D-i$.
  We now prove that the pivot in $\eta_L$ is at index
  $d_y+D-\rho$. Recall that we write $h_1,\dots,h_L$ for the canonical
  lifts of $\eta_1,\dots,\eta_L$ to  $\K[x]^{d_y+D}$; in
  particular, the pivot index $r$ of $h_L$, as defined in the
  proposition, is also the pivot index of $\eta_L$, so our claim is
  $r=d_y+D-\rho$.  
  Suppose that the pivot in $\eta_L$ is at an index different from
  $d_y+D-\rho$. By the previous discussion, it can only lie at a
  larger index, say $m > d_y+D-\rho$; this may happen only if $\rho >
  0$, in which case we saw that
  $\rho=\nu_{\sigma-1}=\deg_y(\Gamma_{\sigma-1})$ and
  $\Gamma_{\sigma}=x^k$.

  Let $H_1,\dots,H_L$ be obtained by applying
  $\pi_{d_y+D}^{-1}$ to $h_1,\dots,h_L$. It follows that $H_L$ has
  $y$-degree $d_y+D-m < \rho=\deg_y(\Gamma_{\sigma-1})$, and $x$-degree
  less than $k = \deg_x(\Gamma_\sigma)$. Thus, $H_L$ is reduced with
  respect to the Gr\"obner basis $\bm\Gamma$ of $J$. On the other
  hand, because $\eta_L$ is in the column span of $\bar{\bm S}$, its
  canonical lift $h_L$ is in the column space of $\bm S$, up to the
  addition of a vector with entries in $x^k \K[x]$. In other words,
  $H_L$ is in $J$, so that $H_L$ must be zero, a contradiction.

  Thus, the pivot index of $\eta_L$ is exactly $d_y+D-\rho$, that is,
  the same as that of $v_{\rho}$. Our previous discussion on the
  pivots in the vectors $\eta_i$ then implies that for
  $i=\rho,\dots,\nu_0$, the pivot index of $\eta_{L+\rho-i}$ is
  $d_y+D-i$, that is, the same as that of $v_i$. This implies that
  \vspace{-0.5em}
  \begin{equation}\label{eq:veta}
  v_i = \sum_{j=\rho}^i \alpha_{i,j} \eta_{L+\rho-j},
  \end{equation}
      for some coefficients $\alpha_{i,j}$ in $\A=\K[x]/x^k$.  On the
  other hand, all polynomials $H_L,\dots,H_{L+\rho-\nu_0}$ are in $J$
  (by the argument we used for $H_L$). By \cref{lemma:freem}, we
  deduce that for $i=\rho,\dots,\nu_0$, $H_{L+\rho-i}$ can be written
  as $H_{L+\rho-i}= \sum_{j=\rho}^i \beta_{i,j} C_j$, for some
  coefficients $\beta_{i,j}$ in $\K[x]$. Applying
  $\pi_{d_y+D}$ and reducing modulo $x^k$, this gives 
  \vspace{-0.5em}
  \begin{equation}\label{eq:etav}
  \eta_{L+\rho-i}= \sum_{j=\rho}^i \bar \beta_{i,j} v_j,     
  \end{equation}
  with $\bar\beta_{i,j}=\beta_{i,j} \bmod x^k$ for all $i,j$.
  We know that the pivots of both $v_i$ and $\eta_{L+\rho-i}$ are
  powers of $x$ (the latter, by the properties of the Howell form), so
  \cref{eq:veta} and \cref{eq:etav} show that the pivots in $v_i$ and
  $\eta_{L+\rho-i}$ are the same, for $i=\rho,\dots,\nu_0$.
  
  Back in $\K[x,y]$, we deduce that $C_i$ and $H_{L+\rho-i}$ have the
  same coefficient in $y^i$, for $i=\rho,\dots,\nu_0$.  As
  in the proof of \cref{detached}, we deduce that we actually have
  $C_i=H_{L+\rho-i}$ for $i=\rho,\dots,\nu_0$: we observe
  that their terms of $y$-degree less than $i$ are reduced with
  respect to $\Gamma$; it follows that $C_i-H_{L+\rho-i}$ is both
  in $J$ and reduced with respect to its lexicographic Gr\"obner
  basis, so it vanishes.

  Taking $i=\nu_0$, we deduce in particular that $H_{L+\rho-\nu_0}$ is
  monic in $y$ (and no $H_i$ of larger index has this property), so
  the index $L'$ defined in the proposition is $L'=L+\rho-\nu_0$; the
  corresponding polynomials are $C_{\nu_0},\dots,C_\rho$.
  Since we saw that $r=d_y+D-\rho$, the integer $L''$ in the
  proposition is $L''=L+\rho$, and through $\pi_{d_y+D}^{-1}$, the
  columns $h_{L+1},\dots,h_{L+\rho}$ become $y^{\rho-1} x^k,\dots,x^k$
  (there is no such column if $\rho=0$). These are precisely the
  polynomials $C_{\rho-1},\dots,C_0$ that were missing if $\rho > 0$.
\end{proof}
\vspace{-0.5em}
We call $\textsc{HowellGroebnerBasis}(\mathcal F, k, D)$ a procedure
that takes as input $\mathcal F=(f_1,\dots,f_t)$, $k$ and $D$, and
returns the lexicographic Gr\"obner basis of $\langle
f_1,\dots,f_t,x^k \rangle$ obtained from the Howell form of $\bar{\bm
  S}$, taking for $d_y$ the maximum of the degrees of $f_1,\dots,f_t$,
and choosing for $D$ the integer prescribed by \cref{DeltaNullst}. In
this case, there is no need to make $\bar{\bm S}$ square: the
algorithm of~\cite[Chapter~4]{Storjohann00} computes its Howell form
using $\softO(tD^\omega k)$ operations in $\K$.

The main application we will make of Howell form computation is to
obtain the Gr\"obner basis of the $\langle x,y\rangle$-primary
component of an ideal such as $I=\langle f_1,\dots,f_t\rangle$. In
order to do so, we will assume that we are in ``nice'' coordinates, in
the sense that the projection on the first factor $V(\mathcal F) \to
\overline \K$ is one-to-one.
\vspace{-0.5em}
\begin{lemma}
  Let $\mathcal F=(f_1,\dots,f_t)$ be in $\K[x,y]$, and suppose that
  the projection on the first factor $V(\mathcal F) \to \overline \K$
  is one-to-one.  Let further $J$ be the $\langle x,y\rangle$-primary
  component of $I=\langle f_1,\dots,f_t\rangle$, with $m$ the smallest
  integer such that $x^m$ is in $J$.  Then: the smallest power of $x$
  in the ideal $H=\langle f_1,\dots,f_t,x^k\rangle$ is
  $x^{\min(m,k)}$, and for $k \ge m$, $H=J$.
\end{lemma}
\vspace{-1em}

\begin{proof}
  First, we establish that $J=\langle
  f_1,\dots,f_t,x^m\rangle$. For one direction, all $f_i$'s, as
  well as $x^m$, are in $J$ by definition.  Conversely, the
  assumption on $V(\mathcal F)$ implies that we can write $\langle
  f_1,\dots,f_t\rangle = J J'$, with $J'$ having no solution above
  $x=0$; in particular, there exist polynomials $u,v$ with $u x^m +
  v = 1$ and $v$ in $J'$. From this, we get $J = (u x^m + v ) J$,
  and every element in $u x^m J$ is a multiple of $x^m$, while
  every element in $vJ$ is in $\langle f_1,\dots,f_t\rangle.$

  Suppose $k \ge m$. We now have polynomials $u',v'$ with $u' x^{k-m}
  + v' = 1$ and $v'$ in $J'$. Multiplying by $x^m$ shows that $x^m$ is
  in the ideal $H=\langle f_1,\dots,f_t,x^k\rangle$, so that $H=J$
  (this proves the last claim in the lemma). In this case, the
  smallest power of $x$ in $H$ is thus $x^m$.

  Suppose $k \le m$. In this case, we prove that the minimal power
  of $x$ in $H=\langle f_1,\dots,f_t,x^k\rangle$ is $x^k$. First, note
  that in this case, $H=\langle
  f_1,\dots,f_t,x^m,x^k\rangle=J+\langle x^k \rangle$, and let
  $x^e$ be the minimum power of $x$ in $H$; suppose $e < k$, so that
  $e < m$.  It follows that $x^e$ is the normal form of a
  polynomial of the form $f x^k$, modulo the Gr\"obner basis $\mathcal G$
  of $J$.  However, Lazard's structure
  theorem~\cite[Theorem~1]{Lazard85} implies that through reduction
  modulo such a Gr\"obner basis, no term of $x$-degree less than $k$
  can appear; a contradiction.  
\end{proof}
\vspace{-0.5em}
This allows us to design an algorithm \textsc{GroebnerBasisAtZero} that
computes the Gr\"obner basis of $J$ (under the position assumption in
the lemma), even though we do not know $m$ in advance: we call
\textsc{HowellGroebnerBasis} with inputs the polynomials
$(f_1,\dots,f_t,x^k)$, for $k=2^i$, with $i=0,1,\dots$, until the
output does {\em not} contain $x^k$. Indeed, the lemma shows that if
$x^k$ is in the Gr\"obner basis of $H=\langle
f_1,\dots,f_t,x^k\rangle$, we have $k \le m$, while if it is not, we
have reached $k > m$, and the output is the Gr\"obner basis of $J$.

Altogether, we do $O(\log(m))$ calls to \textsc{HowellGroebnerBasis},
with $k \le 2m$. With $d$ the maximum degree of $f_1,\dots,f_t$, the
runtime is $\softO(tD^\omega m)$ operations in $\K$, with $D$ in
$\{\Delta_1(d_y),\Delta_2(d),\Delta_3(d)\}$, depending on our
assumptions on $f_1,\dots,f_t$ (recall that $d_y$ and $d$ are the
maximum $y$-degree, resp.\ degree, of the input).

\vspace{-0.5em}
\section{Coefficient size and bad reductions}\label{sec:mod}

Our goal now is to give height bounds on the elements in the
lexicographic Gr\"obner basis of some polynomials $\mathcal
F=(f_1,\dots,f_t)$, working specifically over $\K=\Q$. In this
section, we assume that the input polynomials have integer
coefficients.

The {\em height} $u \in \Z-\{0\}$ is simply $\log(|u|)$. The key
quantity $H(\mathcal F)$, together with a nonzero integer
$\beta_{\mathcal F} \in \Z$, are defined as follows.
\begin{definition}\label{def:H}
  Consider polynomials $\mathcal F= (f_1,\dots,f_t)$ in $\Z[x,y]$, let
  $I$ be the ideal they generate in $\Q[x,y]$, with lexicographic
  Gr\"obner basis $\mathcal G=(g_0,\dots,g_s)$. We define $H(\mathcal
  F)$ as the smallest integer such that there exists $\beta_{\mathcal
    F}$ nonzero in $\Z$ for which we have:
  \begin{itemize}[leftmargin=*]
  \item the polynomials $\beta_{\mathcal F} g_0,\dots,\beta_{\mathcal
    F} g_s$ are in $\Z[x,y]$
  \item all coefficients of $\beta_{\mathcal F}
    g_0,\dots,\beta_{\mathcal F} g_s$ (which include in particular
    $\beta_{\mathcal F}$ itself) have height at most $H(\mathcal F)$
  \item for any prime $p$ in $\Z$, if $p$ does not divide
    $\beta_{\mathcal F}$, $\mathcal G \bmod p$ is the lexicographic
    Gr\"obner basis of $\langle f_1 \bmod p,\dots,f_t \bmod p \rangle$
    in $\F_p[x,y]$.
  \end{itemize}
\end{definition}

In order to give upper bounds on $H(\mathcal F)$, we introduce two
functions $B(n,d,h)$ and $C(t,d,D,h)$. The first one is defined by
$$B(n,d,h)=(N+1) h + N\log(N) + \log(n(d+1)),$$
 with $N= n^2d -nd +n$, whereas $C(t,d,D,h)$ is defined by
$$C(t,d,D,h)=B(tD,d,h)+h + \log(2).$$ In particular, $B(n,d,h)$ is in
 $\softO(n^2dh)$ and $C(t,d,D,h)$ is in $\softO(t^2D^2dh)$.
  
\begin{proposition}\label{prop:heightGB}
  Let $\mathcal F= (f_1,\dots,f_t)$ be in $\Z[x,y]$, for $t \ge 2$,
  such that the ideal $I=\langle f_1,\dots,f_t\rangle \subset \Q[x,y]$
  has dimension zero. Suppose that all $f_i$'s have $y$-degree at most
  $d_y$, degree at most $d$, and coefficients of height at most $h$.
  \begin{enumerate}
\item[{(i)}]   if there exists $i$ in $\{1,\dots,t\}$ such that the
    coefficient of $y^{d_y}$ in $f_i$ is a nonzero constant,
    $H(\mathcal F) \le C(t,d,\Delta_1(d_y),h) \in \softO(t^2d^3h)$
 \item[{(ii)}]   if $t=2$, $H(\mathcal F) \le C(2,d,\Delta_2(d),h) \in
    \softO(d^5h)$
 \item[{(iii)}]    in general, $H(\mathcal F) \le C(t,d,\Delta_3(d),h)
    \in \softO(t^2d^9h)$.
  \end{enumerate}
\end{proposition}
The proposition will follow from height bounds for Hermite forms of
matrices due to Storjohann, with the results in the previous
section. We do not have a direct equivalent for Howell forms (we are
not aware of previous work about height bounds or primes of bad
reduction in this context): if we are interested in the primary
component of $I$ at the origin, we may apply the results above to the
polynomials $f_1,\dots,f_t,x^k$, for a large enough $k$.

To our knowledge, no comparable bounds were given in this
setting. Several previous results discussed the case of {\em radical}
ideal with finitely many solutions. If their Gr\"obner basis $\mathcal
G$ is a {\em triangular set}, the
results in~\cite{DaSc04} show that the polynomials in $\mathcal G$
have coefficients with numerator and denominator of height $\softO(d^3
h+d^4)$. Our result does not feature the term $d^4$, but this might be
due to the proof techniques of~\cite{DaSc04}, which are not limited to
systems in two variables. If we keep the radicality assumption, but
allow arbitrary leading terms, the best previous bound we are aware of
is $\softO(d^7h + d^8)$, from~\cite{Dahan09}.

We start the proof with a result due to Storjohann.
\begin{proposition}[{\cite[Section~6.2]{Storjohann94}}]\label{prop:storj}
  Let ${\bm A}$ be in $\Z[x]^{n\times n}$, with nonzero determinant
  and entries of degree at most $d > 0$ and height at most $h$. Let
  further $\bm H$ be the Hermite normal form of $\bm A$. Then, there
  exists $\alpha$ nonzero in $\Z$ such that all entries of $\alpha \bm
  H$ are in $\Z[x]$, $\alpha$ and the coefficients of all entries of
  $\alpha \bm H$ have height at most $B(n,d,h)$, and for any prime
  $p$, if $p$ does not divide $\alpha$, then $\bm H \bmod p$ is the
  Hermite normal form of $\bm A \bmod p$ in $\F_p[x]^{n\times n}$.
\end{proposition}

Let then $f_1,\dots,f_t$ be as in \cref{prop:heightGB}. First, we
define integers $\gamma$ and $D$ through the following case
discussion. If we are in case $(i)$, we know that at least one of the
  $f_i$'s has a coefficient of $y$-degree $d_y$ in $\Z-\{0\}$; let
  $\gamma$ be such a coefficient. We let $D=\Delta_1(d_y)$ from
  \cref{DeltaNullst}.
In case $(ii)$ or $(iii)$, we let $\gamma=1$, and we take
  respectively $D=\Delta_2(d)$ or $D=\Delta_3(d)$, with notation as above.
In any case, we know that $\Delta(\mathcal F) \le D$, so we can apply
\cref{detached}; it shows that we can recover the (minimal, reduced)
lexicographic Gr\"obner basis of $I=\langle f_1,\dots,f_t\rangle$ from
the columns of the Hermite form of the Sylvester-like matrix $\bm S$
defined in that proposition.

As in the previous section, there is a $(d_y+D) \times tD$ matrix $\bm
S'$ obtained by permuting the columns of $\bm S$ whose leading
$(d_y+D)\times (d_y+D)$ minor is nonzero. Consider again the $tD
\times tD$ square matrix ${\bm S}^{{\rm sq}}$ of \cref{eq:defSsq} and
its Hermite form ${\bm H}^{{\rm sq}}$; the first $d_y+D$ rows of ${\bm
  H}^{{\rm sq}}$ are the Hermite form $\bm H$ of~$\bm S$.

Since ${\bm S}^{{\rm sq}}$ has nonzero determinant, we let $\alpha$ be
the nonzero integer associated to it by \cref{prop:storj}, and set
$\beta = \alpha \gamma$. Then, all entries of $\beta
\bm H^{{\rm sq}}$, and thus of $\beta \bm H$, are in $\Z[x]$, the
latter having coefficients of height at most $C(t,d,D,h)$.  By
\cref{detached},  these bounds apply in
particular to the Gr\"obner basis $(g_0,\dots,g_s)$ of $I$.

Suppose that $p$ is a prime that does not divide $\beta$. Because $p$
does not divide $\alpha$, \cref{prop:storj} shows that $\bar{\bm
  H}^{{\rm sq}}=\bm H^{{\rm sq}} \bmod p$ is the Hermite normal form
of $\bar{\bm S}^{{\rm sq}}=\bm S^{{\rm sq}} \bmod p$. Considering only
the first $tD$ rows, we see that $\bar{\bm H}=\bm H \bmod p$ is the
Hermite normal form of $\bar{\bm S}=\bm S \bmod p$. Now, let us prove
that we still have $\Delta(\bar{\mathcal F}) \le D$.
\vspace{-0.25em}
\begin{itemize}[leftmargin=*]
\item If we are in case $(i)$, since $p$ does not divide $\gamma$, at
  least one of the polynomials $\bar f_i=f_i \bmod p$ has its
  coefficient of $y$-degree $d_y$ a nonzero constant in $\F_p$.  Since
  all $\bar f_i$'s have $y$-degree at most $d_y$, we deduce
  $\Delta(\bar{\mathcal F})=d_y$ in this case (first item of
  \cref{DeltaNullst})
\item If we are in case $(ii)$ or $(iii)$, the discussion above shows
  that $\bar g_0$ and $\bar g_s$ are in the ideal $\langle \bar
  f_1,\dots,\bar f_t\rangle$, so that this ideal admits finitely many
  solutions in an algebraic closure of $\F_p$. Using the second and
  third items of \cref{DeltaNullst} gives our claim.
\end{itemize}
\vspace{-0.25em}
We can then apply \cref{detached} to $\bar {\mathcal F}=(\bar
f_1,\dots,\bar f_t)$, and deduce that the columns of the Hermite form
of $\bar{\bm S}$ give a detaching basis, and in particular the
lexicographic Gr\"obner basis of $\langle \bar f_1,\dots,\bar
f_t\rangle$. This proves the proposition.

\vspace{-0.5em}
\section{Applying changes of coordinates}\label{section:borel-fixed}

Now, we quantify changes of coordinates that ensure desirable
properties. We
write $\bm \gamma$ for a $2\times 2$ matrix $\bm \gamma =
[\gamma_{i,j}]_{1\le i,j\le 2}$ with entries in $\overline\Q$, and we
identity $M_2(\overline\Q)$ with $\overline\Q{}^4$ through $\bm\gamma
\mapsto [\gamma_{1,1},\gamma_{1,2},\gamma_{2,1},\gamma_{2,2}]$.  For
$\gamma$ in ${\rm GL}_2(\overline\Q)$ as above and $f$ in
$\overline\Q[x,y]$, we write $f^{\bm\gamma}=f(\gamma_{1,1} x +
\gamma_{2,1} y, \gamma_{1,2}x + \gamma_{2,2}y)$.

For $\mathcal F=(f_1,\dots,f_t)$ as in the previous sections, the best
degree and height bounds $\Delta(\mathcal F)$ and $H(\mathcal F)$
apply when the input equations have a particular property: at least
one $f_i$ has a term of maximal degree that involves $y$
only. Geometrically, this means that the curve $V(f_i) \subset \overline\Q{}^2$
has no vertical asymptote; we also say that it is in Noether
position. The following lemma is straightforward.
\vspace{-1mm}
\begin{lemma}\label{prop:F1}
  Take $f$ in $\Q[x,y]$ of degree $d$. Then there exists a
  hypersurface $Y_1 \subset \overline\Q{}^4$ of degree at most $d$
  such that if $\bm\gamma$ is in $\overline\Q{}^4-Y_1$, the
  coefficient of $y^d$ in $f^{\bm\gamma}$ is nonzero.
\end{lemma}
\vspace{-0.5em}
Another favorable situation, illustrated when we dealt with Howell
forms, occurs when the projection $V(\mathcal F) \to \overline\Q$
given by $(\alpha,\beta)\mapsto \alpha$ is one-to-one. Again, the
proof is standard (see e.g.~\cite{Rouillier99}), once we see that
$V(\mathcal F)$ has cardinality at most $d^2$.

\vspace{-0.5em}
\begin{lemma}\label{prop:F2}
  Let $\mathcal F= (f_1,\dots,f_t)$ be in $\Q[x,y]$ of degrees at most
  $d$, and suppose that $V(\mathcal F)$ is finite. Then there exists a
  hypersurface $Y_2 \subset \overline\Q{}^4$ of degree at most $d^4$
  such that if $\bm\gamma$ is invertible and in $\overline\Q{}^4-Y_2$,
  the projection on the first factor $V(\bm f^{\bm \gamma}) \to
  \overline\Q$ is one-to-one.
\end{lemma}
\vspace{-0.5em}
\section{Main algorithms}\label{sec:mainalgo}

We can finally present our main algorithms, where we use Newton
iteration to compute lexicographic Gr\"obner bases: we are given
$\mathcal{F}=(f_1,\dots,f_t)$ in $\Z[x,y]$, and we compute either the
Gr\"obner basis $\mathcal G=(g_0,\dots,g_s)$ of $I=\langle
f_1,\dots,f_t \rangle$, or the Gr\"obner basis
$\mathcal{G}^0=(g^0_0,\dots,g^0_r)$ of the $\langle
x,y\rangle$-primary component of $I$ using $p$-adic approximation, for
a prime $p$. In what follows, we give details for the computation of
$\mathcal G$; we will mention what modifications are needed if we want
to compute $\mathcal{G}^0$.

The algorithm is randomized; it takes a parameter $P \ge 1$, our
goal being to obtain the correct output with probability at least
$1-1/2^P$. Throughout, we assume that $f_1$ has maximum degree among
the $f_i$'s (we write $d=\deg(f_1)$) and that $I$ has dimension zero.
Let $\delta =\deg(I) =\dim_\Q \Q[x,y]/I$, $\delta \le d^2$ and let $b$
be the maximum height of the numerators and denominators of the
coefficients in $\mathcal G$. Each polynomial in $\mathcal G$ has at
most $\delta+1$ coefficients, so the total bit-size of the output is
$O(s \delta b)$.

\vspace{-0.5em}
\subsection{Overview}\label{sec:overview}

We start by presenting the main steps of the algorithm, leaving out
some details of the analysis for the next subsection. Runtimes are
given in terms of bit operations; here, we use the fact that
operations $(+,\times)$ modulo a positive integer $M$ take
$\softO(\log(M))$ bit operations, as does inversion modulo $M$ if $M$
is prime~\cite{GaGe13}.

\smallskip\noindent{\bf Introducing a change of coordinates.} We first
choose a change of variables $\bm\gamma$ with coefficients in
$\Z$. Applying it to the input equations $\mathcal F$ gives
polynomials $\mathcal{H}=(h_1,\dots,h_t)$, which we do not need to
compute explicitly (as they may have large height). We let
$\mathcal{B} = (B_0,\dots,B_\sigma)$ be the lexicographic Gr\"obner
basis of these polynomials in $\Q[x,y]$ (as with $\mathcal{H}$, we
do not compute it explicitly).

We assume that $\bm\gamma$ satisfies the assumptions of
\cref{prop:F1,prop:F2}, so that their conclusions hold.

\smallskip\noindent{\bf Computing Gr\"obner bases modulo $\bm p$.}
Next, we choose two primes $p,p'$, and compute the
Gr\"obner bases $\mathcal B_{p}$ of $(\mathcal H \bmod p)$, and $\mathcal B_{p'}$ of
$(\mathcal H \bmod p')$. We assume that neither $p$ nor $p'$ divides
the integers $\beta_{\mathcal F}$ and $\beta_{\mathcal H}$ from
\cref{def:H} applied to $\mathcal F$ and $\mathcal H$, respectively. In
particular, all denominators in $\mathcal B$ are invertible modulo $p$
and $p'$, and $\mathcal B_{p}=\mathcal B \bmod p$ and $\mathcal B_{p'}=\mathcal
H \bmod p'$.

To compute $\mathcal B_{p}$ and $\mathcal B_{p'}$, the algorithm reduces the $O(t d^2)$
coefficients of $\mathcal F$ modulo $p$ and $p'$. Then, we apply
$\bm\gamma$ to the results, to obtain $\mathcal H \bmod p$ and
$\mathcal H \bmod p'$. Due to \cref{prop:F1}, the coefficient of $y^d$
in $h_1$ is a nonzero constant; if this is still the case modulo $p$
and $p'$, we use $\textsc{HermiteGroebnerBasis}$ with $D=d$ to get
$\mathcal B_{p}$ and $\mathcal B_{p'}$; otherwise, we raise an error.

\smallskip\noindent{\bf Changing coordinates in $\mathcal B_{p}$ and $\mathcal B_{p'}$.}
Using the Gr\"obner bases $\mathcal B_{p}$ and $\mathcal B_{p'}$ of $(\mathcal H \bmod
p)$ and $(\mathcal H \bmod p')$, we compute the Gr\"obner bases of
$(\mathcal F \bmod p)$ and $(\mathcal F \bmod p')$. This is done using
the algorithm of~\cite{NeSc20}. Since $pp'$ does not divide
$\beta_{\mathcal F}$, we deduce that we obtain $\mathcal G_{1} = \mathcal G
\bmod p$ and $\mathcal G_{1}' =\mathcal G \bmod p'$.

\smallskip\noindent{\bf Computing $\mathcal G_{k}$.} At each step of the main
loop, we start from $\mathcal G_{{k/2}} = \mathcal G \bmod p^{k/2}$, and we
compute $\mathcal G_{{k}}=\mathcal G \bmod p^{k}$. For this, we first need
$\mathcal F \bmod p^{k}$; then, we use procedure
\textsc{LiftOneStepGroebner} from~\cite[Remark~7.3]{ScSt23} to obtain
$\mathcal G_{{k}}$.

\smallskip\noindent{\bf Rational reconstruction.}  We next attempt to
recover all rational coefficients of $\mathcal G$, given those
of $\mathcal G_{k} = \mathcal G \bmod p^k$. For each coefficient $\alpha$ of
$\mathcal G_k$, we attempt to recover a pair $(\eta,\theta)$ in $\Z\times
\N$, with $|\eta| < p^{k/2}/2$ and $\theta \le p^{k/2}$, $\theta$
invertible modulo $p$ and $\alpha = \eta/\theta \bmod p^k$.

By assumption, all nonzero coefficients of $\mathcal G$
have numerators and denominators of height at most $b$, it follows
that if $p^{k/2} > 2 e^b$, we will succeed and correctly recover the
corresponding coefficient in $\mathcal G$~\cite[Theorem~5.26]{GaGe13}.
For smaller values of $k$, rational reconstruction may find no
solution (in which case we reenter the lifting loop at precision
$2k$), or may already terminate; in this case, its output $\mathcal
G_{\rm rec}$ may be different from $\mathcal G$.


\smallskip\noindent{\bf Testing for correctness.} The final step in
the loop is a randomized test, using $\mathcal G_{1}'=\mathcal G \bmod
p'$ as a witness to detect those cases where rational reconstruction
returned an incorrect result. We attempt to reduce $\mathcal G_{\rm
  rec}$ modulo our second prime $p'$ ($\mathcal G_{\rm red}$); if this fails (because $p'$
divides one of the denominators in it), we reenter the lifting loop at
precision $2k$. We
simply compare $\mathcal G_{\rm red}$ and $\mathcal G_{1}'=\mathcal G \bmod
p'$. If they coincide, we return $\mathcal G_{\rm rec}$, otherwise, we
reenter the lifting loop.
\vspace{-1em}
\begin{algorithm}
 \caption{GroebnerBasis(${\mathcal F}$)}
  \begin{algorithmic}[1]
    \Require $\mathcal F=(f_1,\dots,f_t)$ in $\Z[x,y]$, $d=\max\{\deg f_i\}$, 
    \Ensure the lexicographic Gr\"obner basis of $\mathcal F$ in $\Q[x,y]$
    \State choose $\bm \gamma$ in $M_2(\Z)$
    \State choose two different primes $p,p'$;  do steps 3-6  for $i\in \{p,p'\}$
    \State {\bf if} {$\bm \gamma \bmod i$ is not invertible} {\bf then} {raise an error}
    \State $\mathcal H_i \gets \textsc{ChangeCoordinates}(\mathcal F \bmod i, \bm\gamma \bmod i)$
    \State {\bf if} {the coefficient of $y^d$ in $\mathcal H_i(1)$ is zero} {\bf then} {raise an error}
    \State $\mathcal B_{i} \gets \textsc{HermiteGroebnerBasis}(\mathcal H_{i}, d)$

    \State $\mathcal G_{1} \gets \textsc{ChangeCoordinatesGroebner}(\mathcal B_{p}, \bm \gamma^{-1} \bmod p)$
    \State $\mathcal G_{1}' \gets \textsc{ChangeCoordinatesGroebner}(\mathcal B_{p'}, \bm \gamma^{-1} \bmod p')$

    \State $k \gets 1$
    \Repeat
    \State $k \gets 2k$
    \State $\mathcal G_{k} \gets \textsc{LiftOneStepGroebner}(\mathcal F \bmod p^k, \mathcal G_{{k/2}})$
    \State $error, \mathcal G_{\rm rec} \gets \textsc{RationalReconstruction}(\mathcal G_{k})$
    \State {\bf if\ not} $error$ {\bf then}
    {$error, \mathcal G_{\rm red} \gets \mathcal G_{\rm rec} \bmod p'$}
    \Until \textbf{not} $error$ \textbf{and} $\mathcal G_{\rm red} = \mathcal G_1'$
    \State \Return $\mathcal G_{\rm rec}$
  \end{algorithmic}
\end{algorithm}
\vspace{-1em}\\ 



\subsection{Analysis}\label{sec:analysis}
We assume that choosing a random integer in a set $\{0,\dots,A\}$
(uniform distribution) uses $\softO(\log(A))$ bit operations.
We assume that we have an oracle $\mathcal{O}$, which takes as input an
integer $C$, and returns a prime number in $I=[C+1,\dots,2C]$,
uniformly distributed within the set of primes in $I$, using
$\softO(\log(C))$ bit operations.

\smallskip\noindent{\bf Parameters choice.}  The change of variables
$\bm \gamma$ needs to avoid a hypersurface $Y \subset \overline\Q{}^4$
of degree at most $A_1=d^4+d$. We choose its entries uniformly at
random in $\{0,\dots,2^{P+2} A_1\}$; the cost of getting $\bm\gamma$
is negligible.
Then, by the De Millo-Lipton-Schwartz-Zippel lemma, the probability
that $\bm \gamma$ lies on $Y$ is at most $1/2^{P+2}$. In what follows,
we assume that this is the case, so all polynomials
$\mathcal H=\mathcal F^{\bm\gamma}$ have coefficients of height at
most $h'=h + d (P + 5 + \log(A_1)) \in \softO(h + dP)$.

Let $\beta_{\mathcal F}$ and $\beta_{\mathcal H}$ be the nonzero
integers from \cref{def:H} applied to respectively $\mathcal F$ and
$\mathcal H$, and define
\vspace{-0.5em}
\begin{align*}
C_{\mathcal F}&= C(t,d,\Delta_3(d),h)\in \softO(t^2 d^9 h)\\
C_{\mathcal H}&=C(t,d,\Delta_1(d),h') \in \softO(t^2d^4 h P).  
\end{align*}
\vspace{-1.5em}

\noindent Proposition~\ref{prop:heightGB} proves ${\rm height}(\beta_{\mathcal F}) \le
C_{\mathcal F} $ and ${\rm height}(\beta_{\mathcal H}) \le C_{\mathcal
  H}$. In particular, the height bound $b$ on the coefficients of
$\mathcal G$ satisfies $b \le C_{\mathcal F}$, so $b$ is in
$\softO(t^2 d^9 h)$.
Let $\mu_1$ be the coefficient of $y^d$ in $h_1$, which has height at
most $h'$. Our first requirement on $p$ and $p'$ is that neither of
them divides $\mu=\beta_{\mathcal F}\beta_{\mathcal H}\mu_1$. This is
a nonzero integer, with ${\rm height}(\mu) \le A_2$, where we set $A_2
= C_{\mathcal F} + C_{\mathcal H} + h' \in \softO(t^2 d^9 h P)$.

Finally, we want to ensure that in the verification step, if $\mathcal
G_{\rm rec}$ and $\mathcal G$ differ, their reductions modulo $p'$,
called $\mathcal G_{\rm red}$ and $\mathcal G_{1}'$, differ as well. Below, we
let $k_0$ be the first $k$ which is a power of two and such that, at
step $k$, rational reconstruction correctly computes $\mathcal G_{\rm
  rec}=\mathcal G$. For this, it suffices that $p^{k/2} > 2e^b$, and
one verifies this implies that $k_0 \le 8b \in \softO(t^2 d^9
h)$. Since all indices $k$ we go through are powers of two, there are
at most $\log_2(8b)$ incorrect indices $k$. 



Suppose then that at step $k < k_0$, we have found $\mathcal G_{\rm
  rec}$ with rational coefficients; they all have numerators and
denominators at most $p^{k/2} \le 2 e^b$; on the other hand, the
coefficients of $\mathcal G$ have numerators and denominators at most
$e^b$. If $\mathcal G_{\rm rec}$ and $\mathcal G$ differ, there exists
a monomial whose coefficients in $\mathcal G_{\rm rec}$ and $\mathcal
G$ are different; it suffices that $p'$ does not divide the
numerator of their difference. This number has an absolute value of at most
$4 e^{2b}$.

Taking all $k < k_0$ into account, our last requirement is that $p'$
also not divide a certain nonzero integer $\mu'$ (that depends on
$p$). This integer $\mu'$ has height at most $\log_2(8b)(2b +
\log(4))$, so that ${\rm height}(\mu') \le A_3$, with
$A_3 =\log_2(8C_{\mathcal F})(2C_{\mathcal F} + \log(4)) \in \softO(t^2 d^9
h).$

To summarize, once $\bm\gamma$ avoids $Y$, it suffices that $p$ does
not divide $\mu$ and $p'$ does not divide $\mu\mu'$ to ensure
success. We can then finally make our procedure for choosing $p$ and
$p'$ explicit:
\begin{itemize}[leftmargin=*]
\item Let $B = 2^{P+3}\lceil A_2 \rceil$. We use the oracle $\mathcal{O}$ to obtain
  a uniformly sampled prime number in $[B+1,\dots,2B]$. There are at
  least $B/(2\log(B))$ primes in this interval, and at most
  $\log(\mu)/\log(B)$ of them can divide $\mu$, so the probability
  that $p$ does is at most $2 \log(\mu)/B$, which is at most
  $1/2^{P+2}$.
\item Let $B'= 2^{P+3}\lceil A_2 + A_3\rceil$. We use the oracle $\mathcal{O}$ to
  pick $p'$ in the interval $[B'+1,\dots,2B']$, and as a result, the
  probability that $p'$ divides $\mu \mu'$ is at most $1/2^{P+2}$.
\end{itemize}
Altogether, the probability that $\bm \gamma$ avoids $\mathcal{H}$,
$p$ does not divide $\mu$ and $p'$ does not divide $\mu\mu'$ (and thus
that the algorithm succeeds) is thus at least $1-3/2^{P+2} \ge
1-1/2^P$. 


\smallskip\noindent{\bf Complexity.}  To find $\mathcal B_p$ and
$\mathcal B_{p'}$: reducing the coefficients, changing coordinates and
$\textsc{HermiteGroebnerBasis}$ uses $\softO(td^2(\log(pp'))$,
$\softO(t d^2 (h + \log(pp')))$~\cite[Corollary~9.16]{GaGe13} and
$\softO(t^\omega d^{\omega+1}(\log(pp')))$ bit operations,
respectively. Inverting the $\gamma$ on $\mathcal B_p$ and $\mathcal
B_{p'}$ takes $\softO(\delta^3)$ operations in $\F_{p'}$, which is
$\softO(\delta^3\log(p'))$ bit operations. To compute $\mathcal
G_{k}$: coefficients reduction takes $\softO(td^2 (h + k\log(p) ))$
bit operations. Algorithm \textsc{LiftOneStepGroebner} takes a
one-time cost of $t\delta^\omega \log(p)$ bit operations, plus for
\vspace{-0.25em}
$$\softO (s^2 n_0 m_s + t\delta(d^2 + d m_s + s \delta)) k \log(p))$$ 
\vspace{-0.25em}
bit operations
per iteration. Here, $n_0=\deg_y(g_0)$ and $m_s=\deg_x(g_s)$. Rational reconstruction takes
$\softO(k\log(p))$ bit operations per coefficient, for a total of
$\softO(s \delta k \log(p))$. For the test: reduction modulo $p'$ takes 
$\softO( b + \log(p'))$ bit operations per coefficient, for a total of
$\softO(s \delta(b+ \log(p'))$.

Furthermore, both $\log(p)$ and $\log(p')$ are in $\softO(P +
\log(tdh))$. Besides, the definition of $k_0$ implies that at all
lifting steps, $k \log(p)$ is in $\softO(b + \log(p))$, that is
$\softO(b + P + \log(tdh))$. After some straightforward
simplifications, the runtime becomes softly linear in $t d^2 h$,
$(t^\omega d^{\omega+1} + \delta^\omega)(P + \log(tdh))$ and $(s^2n_0
m_s + t\delta(d^2 + d m_s + s\delta)) (b+P+\log(tdh)) $.

In order to get a better grasp on this runtime, let us assume that $P$
and the number of equations $t$ are fixed constants, and use the upper
bounds $n_0, m_s\le \delta$. This gives a total bound
softly linear in
$$d^2 h + ( d^{\omega+1} + \delta^\omega)\log(h) + (d^2 \delta + d
\delta^2 + s^2\delta^2) (b+\log(h)).$$ The first term is the input
size, the second describes computations done modulo small primes,
and the last one computations are done modulo higher powers of $p$. The output size $O(s \delta b)$ bits.

\vspace{-0.5em}
\subsection{Computing the $\langle x,y\rangle$-primary component}\label{sec:variant}

Finally, we describe how to modify the algorithm if we are only
interested in the Gr\"obner basis $\mathcal G^0=(G^0_0,\dots,G^0_r)$
of the $\langle x,y\rangle$-primary component $J$ of $I$. In what
follows, we let $\eta$ be the degree $J$, and $c$ be the maximum
height of the numerators and denominators of the coefficients of
$\mathcal G^0$ : 
the output total size is $O(r \eta c)$ bits.

As before, we use a change of coordinates $\bm\gamma$, and we call
$\mathcal B^0$ the Gr\"obner basis of $J^{\bm \gamma}$. Then, we use \textsc{GroebnerBasisAtZero} instead of
\textsc{HermiteGroebnerBasis}, modulo $p$ and $p'$. Since we are in
generic coordinates, we can use degree $D=d$, so the runtime is
$\softO(t d^{\omega} m\log(p p'))$ bit operations, where $m$ is the
maximal $x$-degree of the polynomials in $\mathcal B^0$. We will use
the bound $m \le \eta$.
  
Then \textsc{LiftOneStepPunctualGroebnerBasis}
from~\cite[Rk~7.3]{ScSt23} can be used with an initial cost (bit operations) of
$\softO(t \eta^\omega \log(p))$ and $\softO(t \eta^2
k \log(p))$ at the $k^{th}$ iteration. The rest of the algorithm is unchanged except for a slight difference in conditions of success.

Now, $\bm\gamma$ has to avoid a hypersurface $Y'$ of degree at most
$d^4+d$, in order to guarantee that $J^{\bm\gamma}$ satisfies
\cref{prop:F1,prop:F2}.  The primes $p$ and $p'$ should divide the
denominator of no coefficient in $\mathcal G^0$ and $\mathcal B^0$;
besides, these polynomials reduced modulo $p$ (resp.\ $p'$) should
still define the $\langle x,y\rangle$-primary components of $f_1 \bmod
p,\dots,f_t\bmod p$ and $f^{\bm\gamma}_1 \bmod
p,\dots,f^{\bm\gamma}_t\bmod p$ (resp.\ modulo $p'$).

The $\langle x,y\rangle$-primary component of
$\langle f_1,\dots,f_t\rangle$ is the ideal generated by $\mathcal{F}'
= (f_1,\dots,f_t,x^{d^2},y^{d^2})$; similarly for $\mathcal H=
(f_1^{\bm \gamma},\dots,f_t^{\bm \gamma})$, giving us polynomials
$\mathcal H'$. It is then sufficient that neither $p$ nor $p'$ divides
the integers $\beta_{\mathcal F'} \beta_{\mathcal H'}$ from
\cref{def:H}.  Their heights are in $\softO(t^2 d^6h)$ and $\softO(t^2
d^6h')$, where $h'$ is the height bound on $\mathcal H$.

The rest of the analysis is conducted as before. Given a fixed integer
$P$, we deduce that we can compute the Gr\"obner basis $\mathcal G^0$,
with a probability of success of at least $1-1/2^P$, using
$$\softO(t d^2 h + (t d^\omega\eta + \eta^\omega)(P+ \log(tdh)) +  t\eta^2(c+\log(tdh)))$$
bit operations. Assuming $t$ and $P$ are fixed, this is softly
linear in 
$d^2 h + (d^\omega\eta+\eta^\omega)\log(h) + \eta^2 c.$ To wit,
the input size is linear in $dh$ and that the output size is in $O(r
\eta c) \subset O(\eta^2 c)$, with $r$ the number of polynomials in
$\mathcal G^0$ $\eta$ its degree and $c$ the bit-size of its
coefficients.\\

\smallskip\noindent{\bf  Acknowledgments.}  We thank Arne Storjohann and Vincent Neiger for answering our
questions on Hermite normal form computations. Schost is supported by
an NSERC Discovery Grant. St-Pierre thanks NSERC, Alexander Graham
Bell Canada Graduate Scholarship, FQRNT and the European Research
Council (ERC) under the European Union’s Horizon Europe research and
innovation programme, grant agreement 101040794 (10000 DIGITS) for
their generous support.

\balance

\bibliographystyle{plain}
\bibliography{biblio}

\end{document}